\documentclass{amsart}
\usepackage{amssymb, amsthm, amsmath}

\newtheorem{thm}{Theorem}[section]
\newtheorem{lem}[thm]{Lemma}
\newtheorem{prop}[thm]{Proposition}

\theoremstyle{definition}
\newtheorem{defn}[thm]{Definition}

\theoremstyle{remark}
\newtheorem{rem}[thm]{Remark}

\def\CM{Cohen-Macaulay }

\def\Proj{\mathop{\mathrm{Proj}}}
\def\Spec{\mathop{\mathrm{Spec}}}
\def\Cl{\mathop{\mathrm{Cl}}}
\def\zero{^{\circ}}
\def\mat{{\mathfrak{M}}_{m,n}}

\begin{document}

\title{F-regularity does not deform}

\author{Anurag K. Singh}

\address{Department Of Mathematics, The University Of Michigan, East Hall, 525
East University Avenue, Ann Arbor, MI 48109-1109 
\newline
{\it Present address:} Department of Mathematics, University of Illinois, 1409
W. Green Street Urbana, IL 61801}

\email{singh6@math.uiuc.edu}

\thanks{Manuscript received March 24, 1998} 

\begin{abstract}  
We show that the property of F-regularity does not deform, and thereby settle a
longstanding open question in the theory of tight closure. Specifically, we
construct a three dimensional ${\mathbb N}$-graded domain $R$ which is not 
F-regular (or even F-pure), but has a quotient $R/tR$ which is F-regular.
Examples are constructed over fields of characteristic $p>0$, as well as over
fields of characteristic zero. 
\end{abstract}

\maketitle

\section{Introduction}

Throughout this paper, all rings are commutative, Noetherian, and have an 
identity element. The theory of {\it tight closure}\/ was developed by Melvin
Hochster and Craig Huneke in \cite{HHjams} and draws attention to rings which
have the property that all their ideals are tightly closed, called {\it weakly
F-regular} rings. The term {\it F-regular}\/ is reserved for rings all of 
whose localizations are weakly F-regular. A natural question that arose with
the development of the theory was whether the property of F-regularity deforms,
i.e., if $(R,m,K)$ is a local ring such that $R/tR$ is F-regular for some 
nonzerodivisor $t\in m$, must $R$ be F-regular? (See the Epilogue of
\cite{Ho}.) Hochster and Huneke showed that this is indeed true if the ring $R$
is Gorenstein, \cite{HHbasec}, and their work has been followed by various
attempts at extending this result, see \cite{AKM, Si, Smdeform}. Our primary
goal here is to settle this question by constructing a family of examples to
show that F-regularity does not deform. We shall throughout be considering
$\mathbb N$-graded rings, but local examples can be obtained, in all cases, by
localizing at the homogeneous maximal ideals. Our main result is:

\begin{thm}
There exists an $\mathbb N$-graded ring $R$ of dimension three (finitely
generated over a field $R_0 = K$ of characteristic $p > 2$) which is not 
F-pure, but has an F-regular quotient $R/tR$ where $t\in m$ is a homogeneous
nonzerodivisor.

Specifically, for positive integers $m$ and $n$ satisfying $m - m/n > 2$, 
consider the ring 
$R=K[A,B,C,D,T]/I$ where $I$ is generated by the size two minors of the matrix 
$$ 
\begin{pmatrix} 
A^2+T^m \ & B \ & D \\
C \ & A^2 \ & B^n-D \\ 
\end{pmatrix}. 
$$ 
Then the ring $R/tR$ is F-regular, whereas $R$ is not F-regular. If $p$ and $m$
are relatively prime integers, then the ring $R$ is not F-pure. 
\label{main}
\end{thm}

A notion closely related to (and frequently the same as) F-regularity is that
of {\it strong F-regularity}. A recent result of G.~Lyubeznik and K.~E.~Smith
states that the properties of weak F-regularity, F-regularity, and strong
F-regularity agree for $\mathbb N$-graded F-finite rings, see \cite{LS}. In the
light of this result, we frequently make no distinction between these notions.

The first result on the deformation of F-regularity is a theorem of Hochster
and Huneke which states that for a Gorenstein local ring $(R,m,K)$, if $R/tR$
is F-regular for some nonzerodivisor $t\in m$, then $R$ is F-regular. They show
that for Gorenstein rings the properties of F-regularity and F-rationality
coincide, and that F-rationality deforms, see \cite[Theorem 4.2]{HHbasec}. This
result is generalized in \cite{Si} where the author uses the idea of passing
to an anti-canonical cover $S=\oplus_{i \ge 0}I^{(i)}$, where $I$ represents
the inverse of the canonical module in $\Cl(R)$. Strong F-regularity is shown
to deform in the case that the symbolic powers $I^{(i)}$ satisfy the Serre
condition $S_3$ for all $i \ge 0$, and the ring $S$ is Noetherian.

This question has also been settled for $\mathbb Q$-{\it Gorenstein} rings,
i.e., rings in which the canonical module is a torsion element of the divisor
class group. For $\mathbb Q$-Gorenstein rings essentially of finite type over a
field of characteristic zero, K.~E.~Smith showed that the property of F-regular
type (a characteristic zero analogue of F-regularity) does deform, see
\cite{Smdeform}. The point is that in this setting F-regular type is equivalent
to log-terminal singularities, and log-terminal singularities deform by
J.~Koll\'ar's result on \lq\lq inversion of adjunction\rq\rq, see \cite{Kollar}.
For $\mathbb Q$-Gorenstein rings of characteristic $p$, a purely algebraic
proof that F-regularity deforms was provided by I.~Aberbach, M.~Katzman, and
B.~MacCrimmon in \cite{AKM}.

Before proceeding with formal definitions in the next section, we would like to
point out that although tight closure is primarily a characteristic $p$ notion,
it has strong connections with the study of singularities of algebraic
varieties over fields of characteristic zero. Specifically, let $R$ be a ring
which is essentially of finite type over a field of characteristic zero; then
$R$ has rational singularities if and only if it is of F-rational type, see
\cite{Hara, Smratsing}. In the $\mathbb Q$-Gorenstein case, we have some even
more remarkable connections: F-regular type is equivalent to log-terminal
singularities and F-pure type implies (and is conjectured to be equivalent to)
log-canonical singularities, see \cite{Smvanish, Walog}.

\section{Frobenius closure and tight closure}

By an ${\mathbb N}$-{\it graded ring} $R$, we will always mean a ring
$R=\oplus_{n \ge 0}R_n$, finitely generated over a field $R_0=K$.

Let $R$ be a Noetherian ring of characteristic $p > 0$. The letter $e$ denotes 
a variable nonnegative integer, and $q$ its $e$th power, i.e., $q=p^e$. For an
ideal $I=(x_1,\dots, x_n) \subseteq R$, let $I^{[q]} = (x_1^q,\dots, x_n^q)$.

For a reduced ring $R$ of characteristic $p > 0$, $R^{1/q}$ shall denote the
ring obtained by adjoining all $q$th roots of elements of $R$. A ring $R$ is
said to be {\it F-finite}\/ if $R^{1/p}$ is module-finite over $R$. A finitely
generated algebra $R$ over a field $K$ is F-finite if and only if $K^{1/p}$ is
a finite field extension of $K$. We use $R\zero$ to denote the complement of
the union of the minimal primes of $R$.

\begin{defn} 
Let $R$ be a ring of characteristic $p$, and $I$ an ideal of $R$. For an
element $x$ of $R$, we say that $x \in I^F$, the {\it Frobenius closure }\/ of 
$I$, if there exists $q=p^e$ such that $x^q \in I^{[q]}$. 

An element $x$ of $R$ is said to be in $I^*$, the {\it tight closure}\/ of $I$,
if there exists $c \in R\zero$ such that $cx^q \in I^{[q]}$ for all $q=p^e \gg
0$. If $I=I^*$ we say that the ideal $I$ is {\it tightly closed}. It is easily
seen that $I \subseteq I^F \subseteq I^*$.

A ring $R$ is said to be {\it F-pure }\/ if the Frobenius homomorphism is pure,
i.e., $F: M \to F(M)$ is injective for all $R$-modules $M$. Note that this
implies $I^F = I$ for all ideals $I$ of $R$.

A ring $R$ is {\it weakly F-regular}\/ if every ideal of $R$ is tightly closed,
and is {\it F-regular}\/ if every localization is weakly F-regular. An F-finite
ring $R$ is {\it strongly F-regular}\/ if for every element $c \in R\zero$,
there exists an integer $q=p^e$ such that the $R$-linear inclusion $R \to
R^{1/q}$ sending $1$ to $c^{1/q}$ splits as a map of $R$-modules. $R$ is {\it
F-rational}\/ if, in every local ring of $R$, all ideals generated by systems
of parameters are tightly closed. 
\end{defn}

It follows from the definitions that a weakly F-regular ring is F-rational as
well as F-pure. We summarize some basic results regarding these notions from
\cite[Theorem 3.1]{HHstrong}, \cite[Theorem 4.2]{HHbasec} and 
\cite[Corollaries 4.3, 4.4]{LS}.

\begin{thm} 
\item $(1)$\quad Regular rings are F-regular; if they are F-finite, they are
also strongly F-regular. Strongly F-regular rings are F-regular.

\item $(2)$\quad Direct summands of F-regular rings are F-regular.

\item $(3)$\quad F-rational rings are normal. An F-rational ring which is a
homomorphic image of a \CM ring is itself Cohen-Macaulay.

\item $(4)$\quad An F-rational Gorenstein ring is F-regular. If it is 
F-finite, then it is also strongly F-regular. 

\item $(5)$\quad The notions of weak F-regularity and F-regularity agree for 
${\mathbb N}$-graded rings. For F-finite $\mathbb N$-graded rings, these are
also equivalent to strong F-regularity.
\label{longlist}
\end{thm}

\section{A review of rational coefficient Weil divisors}

The examples constructed in the following section are best understood in the
setting of ${\mathbb Q}$-divisors. Also, an interpretation of the graded pieces
of certain local cohomology modules using ${\mathbb Q}$-divisors provided the
original heuristic ideas which led to these examples. We recall some notation
and results from \cite{De, Wadem, Wadim2}.

\begin{defn} 
By a {\it rational coefficient Weil divisor}\/ (or a ${\mathbb Q}$-{\it
divisor}) on a normal projective variety $X$, we mean a linear combination of
codimension one irreducible subvarieties of $X$, with coefficients in ${\mathbb
Q}$. For $E = \sum n_iV_i$ with $n_i \in {\mathbb Q}$, we set $[E]= \sum
[n_i]V_i$, where $[n]$ denotes the greatest integer less than or equal to $n$,
and define ${\mathcal{O}}_X(E) ={\mathcal{O}}_X([E])$.

Let $E=\sum(p_i/q_i)V_i$ where the integers $p_i$ and $q_i$ are relatively
prime and $q_i > 0$. We define $E' = \sum((q_i-1)/q_i)V_i$ to be the {\it
fractional part\/} of $E$. Note that with this definition we have $-[-nE] =
[nE+E']$ for any integer $n$. 
\end{defn}

For an ample ${\mathbb Q}$-divisor $E$ (i.e., $NE$ is an ample Cartier divisor
for some $N \in {\mathbb N}$), we construct the {\it generalized section ring}:
$$
S = S(X,E) = \oplus _{n \ge 0} H^0(X,{\mathcal {O}}_X(nE)).
$$
In this notation, Demazure's result (\cite[3.5]{De}) states that every 
${\mathbb N}$-graded normal ring arises as a generalized section ring $S =
S(X,E)$ where $E$ is an ample ${\mathbb Q}$-divisor on $X = \Proj S$.

Let $X$ be a smooth projective variety of dimension $d$ with canonical divisor
$K_X$, and let $E$ be an ample ${\mathbb Q}$-divisor on $X$. If $\omega$ denotes 
the graded canonical module of the ring $S=S(X,E)$, K.-i.~Watanabe showed in 
\cite{Wadem} and \cite{Wadim2} that
\begin{gather*}
[\omega^{(i)}]_n = H^0(X,{\mathcal O}_X(i(K_X+E')+nE)) \\
\intertext{and}
[H^{d+1}_m(\omega^{(i)})]_n = H^d(X,{\mathcal O}_X(i(K_X+E')+nE)). 
\end{gather*}
Let $E_S(K) = H^{d+1}_m(\omega)$ denote the injective hull of $K$ as a graded 
$S$-module. The Frobenius action on the $n$th graded piece of $E_S(K)$ can then
be identified with 
$$
H^d(X,{\mathcal O}_X(K_X+E'+nE)) \overset{F}\longrightarrow 
H^d(X,{\mathcal O}_X(p(K_X+E'+nE))), 
$$
and in particular the Frobenius action on $[H^{d+1}_m(\omega)]_0$, the socle of 
$E_S(K)$, can be identified with 
$$
H^d(X,{\mathcal O}_X(K_X+E')) \overset{F}\longrightarrow 
H^d(X,{\mathcal O}_X(p(K_X+E'))). 
$$
If the ring $S$ is F-pure this Frobenius action must be injective, and
consequently $H^d(X,{\mathcal O}_X(p(K_X+E')))$ must be nonzero. Heuristically,
the dimension of the vector space $H^d(X,{\mathcal O}_X(p(K_X+E')))$ may be
regarded as a measure of the F-purity of $S$. With this is mind, we choose $E$
such that $h^d (X,{\mathcal O}_X(p(K_X+E')))$ while nonzero, is \lq\lq
small\rq\rq. This motivates our choice of the ${\mathbb Q}$-divisor $E$ on
${\mathbb P}^1$, see the proof of Proposition \ref{f-reg}.

\section{The main construction}
When working with quotients of polynomial rings, we shall use lower-case letters 
to denote the images of the corresponding variables. 

\begin{rem}
Let $K$ be a field of characteristic $p$. For positive integers $m$ and $n$, 
consider the ring $R=K[A,B,C,D,T]/I$ where $I$ is generated by the size two 
minors of the matrix 
\begin{equation*} 
\mat = \begin{pmatrix} 
A^2+T^m \ & B \ & D \\
C \ & A^2 \ & B^n-D \\ 
\end{pmatrix}. 
\end{equation*}
The ring $R$ is graded by setting the weights of $a$, $b$, $c$, $d$, and $t$ to
be $m$, $2m$, $2m$, $2mn$, and $2$ respectively. This ring is the specialization 
of a \CM ring, and so is itself Cohen-Macaulay. The elements $t$, $c$ and $d$ 
form a homogeneous system of parameters for $R$, and so the element $t\in m$ is 
indeed a nonzerodivisor.
\label{mat}
\end{rem}

We record the following crucial lemma.

\begin{lem} 
Let $m$ and $n$ be positive integers satisfying $m - m/n > 2$. Consider the 
ring $R=K[A,B,C,D,T]/I$ where $I$ is generated by the size two minors of the
matrix $\mat$ (see \S \ref{mat}). If $k$ is a positive integer such that 
$k(m - m/n - 2) \ge 1$, then
$$
(b^nt^{m-1})^{2mk+1} \in (a^{2mk+1}, \ d^{2mk+1}).
$$
\label{keylemma}
\end{lem}

\begin{proof}
Let $\tau=A^2+T^m$ and $\alpha=A^2$. It suffices to working in the polynomial 
ring $K[\tau, \alpha, B, C, D]$ and establish that
$$
B^{n(2mk+1)}(\tau-\alpha)^{2k(m-1)} \in (\alpha^{mk+1}, \ D^{2mk+1})+\mathfrak{a}
$$
where $\mathfrak{a}$ is the ideal generated by the size two minors of the matrix
$$
\begin{pmatrix} 
\tau \ & B \ & D \\
C \ & \alpha \ & B^n-D \\
\end{pmatrix}.
$$
Taking the binomial expansion of $(\tau-\alpha)^{2k(m-1)}$, it suffices to show 
that 
$$
B^{n(2mk+1)}(\tau,\alpha)^{2k(m-1)} \in (\alpha^{mk+1},\ D^{2mk+1})+\mathfrak{a}.
$$
This would follow if we show that for all integers $i$ where $1 \le i \le mk+1$, 
we have
$$
B^{n(2mk+1)}\alpha^{mk+1-i}\tau^{mk-2k+i-1} \in (\alpha^{mk+1}, 
\ D^{2mk+1})+\mathfrak{a},
$$
and so it is certainly enough to show that 
$$
B^{n(2mk+1)}\tau^{mk-2k+i-1} \in (\alpha^{i}, \ D^{2mk+1}) + \mathfrak{a}. 
$$ 
Since $\alpha D-B(B^n-D) \in \mathfrak{a}$, it suffices to establish that
$$
B^{n(2mk+1)}\tau^{mk-2k+i-1} \in (B^{i}(B^n-D)^{i}, \ D^{2mk+1},
\ B^n\tau-D(C+\tau)).
$$
Working modulo the element $B^{i}(B^n-D)^{i}$, we may reduce $B^{n(2mk+1)}$
to a polynomial in $B$ and $D$ such that the highest power of $B$ that occurs is 
less than $i(n+1)$. Consequently it suffices to show that
$$
B^{n(2mk+1-j)}\tau^{mk-2k+i-1}D^{j} \in (D^{2mk+1}, \ B^n\tau-D(C+\tau))
$$
where $n(2mk+1-j) < i(n+1)$, i.e., $j \ge 2mk + (1-i)(1+1/n)$. With this 
simplification, it is enough to check that
$$
B^{n(2mk+1-j)}\tau^{mk-2k+i-1} \in (D^{2mk+1-j}, \ B^n\tau-D(C+\tau)).
$$
It only needs to be verified that $mk-2k+i-1 \ge 2mk+1-j$ since, working modulo 
$B^n\tau-D(C+\tau)$, we can then express $B^{n(2mk+1-j)}\tau^{mk-2k+i-1}$ as a 
multiple of $D^{2mk+1-j}$. Finally, note that
$$
(mk-2k+i-1)-(2mk+1-j) = j-mk-2k+i-2 \ge k(m - \frac{m}{n} - 2) - 1 \ge 0
$$
since $i \le mk+1$, $j \ge 2mk + (1-i)(1+1/n)$ and $k(m - m/n - 2) \ge 1$.
\end{proof}

\begin{prop}
Let $S = K[A,B,C,D]/J$ where the characteristic of the field $K$ is a prime
$p > 2$, and $J$ is the ideal generated by the size two minors of the matrix 
$$
\begin{pmatrix} 
A^2 \ & B \ & D \\
C \ & A^2 \ & B^n-D \\
\end{pmatrix}.
$$
Then $S$ is an F-regular ring. 
\label{f-reg}
\end{prop}

\begin{proof}
There are various ways to establish this. We can identify $S$ with the 
generalized section ring $\oplus_{i \ge 0} H^0(\mathbb P^1,\mathcal O_{\mathbb 
P^1}(iE))X^i$, where $\mathbb P^1 = \Proj K[X,Y]$, and $E$ is the rational 
coefficient Weil divisor 
$$
E = \frac{1}{2}V(X)+\frac{1}{2}V(Y)+\frac{1}{2n}V(X+Y). 
$$
Under this identification, 
$$ 
A = X, \ \ B = \frac{X^3}{Y}, \ \ C = XY \ \text{ and } 
D = \frac{X^{3n+1}}{Y^n(X+Y)}. 
$$ 
One may now appeal to Watanabe's classification in \cite{Wadim2} to conclude
that $S$ is F-regular.

For an alternate proof, it is easily verified that $S$ is the Veronese subring
$$
H^{(2n+1)}=\bigoplus_{i \in \mathbb N}[H]_{i(2n+1)}
$$ 
where $H$ is the hypersurface 
$$ 
K[A, X, Y]/(A^2 - XY(X^n - Y))
$$ 
and the variables $A$, $X$ and $Y$ have weights $2n+1$, $2$ and $2n$
respectively. Here $B = XY^2$, $C = X(X^n - Y)^2$ and $D = Y^{2n+1}$. Since the
characteristic of $K$ is greater than $2$, a routine computation shows that the
hypersurface $H$ is F-regular, and consequently its direct summand $S$ is also
F-regular. 
\end{proof}

\begin{rem}
Although we do not use this fact, we mention that the hypersurface 
$H$ in the proof above is obtained as the cyclic cover 
$$
S \oplus \omega \oplus \omega^{(2)} \oplus \dots \oplus \omega^{(2n)}
$$
where $\omega$ is the canonical module of the ring $S$.    
\end{rem}

\begin{prop}
Let $K$ be a field of characteristic $p > 2$ and consider the ring $R = R_{m,n}= 
K[A,B,C,D,T]/I$ where $I$ is generated by the size two minors of the matrix 
$\mat$ (see \S \ref{mat}). If $m - m/n > 2$, then $R$ is not F-regular. If in 
addition $p$ and $m$ are relatively prime, then $R$ is not F-pure.
\label{non-f-reg} 
\end{prop}

\begin{proof}
First note that $b^nt^{m-1} \notin (a, d)$. To establish that $R$ is not 
F-regular we shall show that $b^nt^{m-1} \in (a, d)^*$. 

For a suitably large arbitrary positive integer $e$, let $q=p^e=2mk+\delta$ 
where $k$ and $\delta$ are integers such that $k(m - m/n - 2) \ge 1$ and $-2m+2
\le \delta \le 1$. To see that $b^nt^{m-1} \in (a, d)^*$, it suffices to show
that 
$$
(b^nt^{m-1})^{q+2m-1} \in (a^q, d^q). 
$$
Since $q+2m-1= 2mk+\delta+2m-1 \ge 2mk+1$ and $q \le 2mk+1$, it suffices to 
check that
$$
(b^nt^{m-1})^{2mk+1} \in (a^{2mk+1}, \ d^{2mk+1}),
$$
but this is precisely the assertion of Lemma \ref{keylemma}.

For the second assertion, note that since $p > 2$, the integers $p$ and $2m$ 
are relatively prime and we may choose a positive integer $e$ such that
$q=p^e=2mk+1$ for some $k > 0$. Taking a higher power of $p$, if necessary, we
may also assume that $k(m - m/n - 2) \ge 1$. But now $(b^nt^{m-1})^{q} \in
(a^q, d^q)$ by Lemma \ref{keylemma}, and so $b^nt^{m-1} \in (a, d)^F$. Hence the
ring $R$ is not F-pure. 
\end{proof}

\begin{proof}[Proof of Theorem \ref{main}]
We have already noted in \S \ref{mat} that the element $t$ is a nonzerodivisor
in $R$, and Proposition \ref{f-reg} establishes that the ring $R/tR$ is 
F-regular. Since $m - m/n >2$, Proposition \ref{non-f-reg} shows that $R$ fails 
to be F-regular, and is not even F-pure if $p$ and $m$ are relatively prime
integers. 
\end{proof}

\section{The characteristic zero case}

Hochster and Huneke have developed a notion of tight closure for rings
essentially of finite type over fields of characteristic zero, see 
\cite{HHjams, HHchar0}. However we can also define notions corresponding to
F-regularity, F-purity, and F-rationality in characteristic zero, without 
explicitly considering a closure operation for rings of characteristic zero. We
include a brief summary, and discuss how the examples constructed above also
show that the property {\it F-regular type}\/ does not deform.

Suppose $R=K[X_1, \dots, X_n]/I$ is a ring finitely generated over a field $K$
of characteristic zero, choose a finitely generated $\mathbb Z$-algebra $A$
such that 
$$
R_A = A[X_1, \dots, X_n]/I_A
$$ 
is a free $A$-algebra with $R \cong R_A \otimes_A K$. Note that the fibers of
the homomorphism $A \to R_A$ over maximal ideals of $A$ are finitely generated
algebras over fields of positive characteristic. 

\begin{defn}
Let $R$ be a ring which is finitely generated over a field of characteristic 
zero. Then $R$ is said to be of {\it F-regular type}\/ if there exists a 
finitely generated $\mathbb Z$-algebra $A \subseteq K$ and a finitely generated
$A$-algebra $R_A$ as above such that $R \cong R_A \otimes_A K$ and, for all
maximal ideals $\mu$ in a Zariski dense subset of $\Spec A$, the fiber rings
$R_A \otimes_A A/\mu$ are F-regular.

Similarly, $R$ is said to be of {\it F-pure type}\/ if for all maximal ideals 
$\mu$ in a Zariski dense subset of $\Spec A$, the fiber rings $R_A \otimes_A
A/\mu$ are F-pure.
\end{defn}

\begin{rem}
Some authors use the term F-pure type (F-regular type) to mean that 
$R_A \otimes_A A/\mu$ is F-pure (F-regular) for all maximal ideals $\mu$ in a 
Zariski dense {\it open}\/ subset of $\Spec A$. 
\end{rem}

\begin{thm}
For positive integers $m$ and $n$ satisfying $m - m/n > 2$, consider the ring 
$R={\mathbb Q}[A,B,C,D,T]/I$ where $I$ is generated by the size two minors of 
the matrix $\mat$ of \S \ref{mat}. Then the ring $R$ is not of F-pure type, 
whereas $R/tR$ is of F-regular type.
\end{thm}

\begin{proof}
If $p$ is a prime integer which does not divide $2m$, the fiber of $\mathbb Z
\to R_{\mathbb Z}$ over $p{\mathbb Z}$ is not F-pure by Proposition 
\ref{non-f-reg}, and consequently the ring $R$ is not of F-pure type. On the
other hand, Proposition \ref{f-reg} shows that $R/tR$ is of F-regular type
since the fiber of ${\mathbb Z} \to (R/tR)_{\mathbb Z}$ over $p{\mathbb Z}$ is
F-regular for all primes $p > 2$.  
\end{proof}

\begin{rem}
R.~Fedder first constructed examples to show that F-purity does not deform, 
see \cite{Fed-pure}. However Fedder pointed out that his examples were less
than satisfactory in two ways: firstly the rings were not integral domains, and
secondly his arguments did not work in the characteristic zero setting, i.e.,
did not comment on the deformation of the property {\em F-pure type}. In
\cite{Si} the author constructed various examples which overcame both these
shortcomings, but left at least one issue unresolved --- although the rings $R$
were integral domains (which were not F-pure), the F-pure quotient rings $R/tR$
were not integral domains. The examples we have constructed here also settle
this remaining issue. 
\end{rem}

\section{Conditions on fibers}

The examples constructed in the previous section are also relevant from the
point of view of the behavior of F-regularity under base change. We first
recall a theorem of Hochster and Huneke, \cite[Theorem 7.24]{HHbasec}.

\begin{thm}
Let $(A,m,K) \to (R,n,L)$ be a flat local homomorphism of local rings of
characteristic $p$ such that $A$ is weakly F-regular, $R$ is excellent, and
the generic and closed fibers are regular. Then the ring $R$ is weakly 
F-regular.
\end{thm}

It is a natural question to ask what properties are inherited by an excellent
ring $R$ if, as above, $(A,m,K) \to (R,n,L)$ is a flat local homomorphism, the
ring $A$ is F-regular and the generic and closed fibers are F-regular. Our 
examples can be used to show that even if $(A,m,K)$ is a discrete valuation
ring and the generic and closed fibers of $(A,m,K) \to (R,n,L)$ are F-regular, 
then the ring $R$ need not be F-regular. 

Once again, we construct $\mathbb N$-graded examples, and examples with local
rings can be obtained by the obvious localizations at the homogeneous maximal
ideals. Let $A = K[T]$ be a polynomial ring in one variable, and 
$R=K[A,B,C,D,T]/I$ where $I$ is generated by the size two minors of the matrix 
$\mat$, see \S \ref{mat}. As before, $K$ is a field of characteristic $p > 2$,
and $m$ and $n$ are positive integers such that $m - m/n > 2$.

The generic fiber of the inclusion $A \to R$ is a localization of $R_t$,
whereas the fiber over the homogeneous maximal ideal of $A$ is $R/tR$. We have
earlier established that $R/tR$ is F-regular, and only need to show that the
ring $R_t$ is F-regular. In the following proposition we show that the $R$ is,
in fact, locally F-regular on the punctured spectrum.

\begin{prop}
Let $K$ be a field of characteristic $p > 2$. For positive integers $m$ and $n$ 
consider the ring $R = R_{m,n}= K[A,B,C,D,T]/I$ where $I$ is generated by the 
size two minors of the matrix $\mat$ (see \S \ref{mat}). Then the ring $R_P$ is 
F-regular for all prime ideals $P$ in $\Spec R - \{m\}$.
\end{prop}

\begin{proof}
A routine verification shows that the singular locus of $R$ is $V(J)$ where 
the defining ideal is $J = (a, \ b, \ c(c+t^m), \ d)$. Consequently we need to 
show that the two local rings $R_P$ and $R_Q$ are F-regular where 
$P=(a,\ b, \ c, \ d)$ and $Q=(a, \ b, \ c+t^m, \ d)$.

After localizing at the prime $P$, we may write $d=b^n(a^2+t^m)/(c+a^2+t^m)$ 
and so $R_P$ is a localization of the ring
$$
K[T, A, B, C]/(A^2(A^2+T^m) - BC)
$$ 
at the prime ideal $(a, \ b, \ c)$. Since $a^2+t^m$ is a unit, the hypersurface
$R_P$ is easily seen to be F-regular. 

Localizing at the prime $Q$, we have $b=a^2(a^2+t^m)/c$ and so  $R_Q$ is a
localization of the ring
$$
K[T, A, C, D]/( C^n D(C+A^2+T^m) - A^{2n} (A^2+T^m)^{n+1} )
$$ 
at the prime ideal $(a, \ c+t^m, \ d)$. Again we have a hypersurface which, it
can be easily verified, is F-regular.  
\end{proof}

\section*{Acknowledgments}

It is a pleasure to thank Melvin Hochster for several enjoyable discussions on
tight closure theory. I am also grateful to Ian Aberbach for pointing out the
issues that lead to \S 6. of this paper.

\end{document}